\documentclass[12pt]{article}
\newcommand{\di}{\displaystyle}

\usepackage[dvips]{graphicx}
\usepackage{verbatim}
\usepackage{amsmath}
\begin{document}
\begin{center}
{\Large\bf Polynomial splines interpolating prime series}\\[0.5cm]

{\large
 L. Alexandrov\footnote[1]{Joint Inst for
 Nuclear Res, BLTP, Dubna, Russia (alexandr@thsun1.jinr.ru,
baranov@thsun1.jinr.ru)}
 D. B. Baranov$^*$
and
 P. T. Yotov\footnote[2]{Institute for Nuclear Research and Nuclear
 Energy, BAS, Sofia (pyotov@inrne.bas.bg)}}
\end{center}

\begin{quote}
{\bf Abstract:}
Differentiable real function
reproducing primes up to a given number and
having a differentiable inverse function is constructed.
This inverse function is compared with the Riemann-Von Mangoldt
exact expression for the number of primes not exceeding
a given value.
Software  for computation of the direct
and inverse functions and their derivatives is developed.
Examples of approximate solution of Diophantine equations
on the primes are given.
\end{quote}

 {\bf 1. Introduction}\\

This article introduces real functions reproducing the values
of mutually inverse arithmetic functions $p(n):\, N\to P$ ({\it
prime $p(n)$ at a number $n$}) and $p^{-1}(p):\, P\to N$ ({\it
number $n$ of the prime $p(n)$}).

The found functions are employed to create
subroutines  for computation of $p(x),\, \frac{\di dp(x)}{\di dx}$ on
$1\leq x<\infty$, and $p^{-1}(x)$, $\frac{\di dp^{-1}(x)}{\di dx}$
on $2\leq x< \infty$.\\

The above-noted programs can be used for a numerical solution of different
problems on the set of primes $P$, including
approximate solution of Diophantine equations on $P$.

The idea consists in establishing a differentiable function which
would include the values of primes and which would
allow one to construct  an inverse function $p^{-1}(x)$
by a Newton method.
More precisely,   the sought function   $p(x),\, 1\leq x< \infty$
should satisfy the following conditions:

a) $p(x)$ reproduces primes;

b) there exists a positive derivative  $\frac{\di dp(x)}{\di dx}$;

c) there exists an inverse function  $p^{-1}(x):\, [2,\infty)\to
[1,\infty)$.

As is known, there are no one-variable polynomials which can produce all the
primes, or primes only. However, this article shows
that there exist polynomial splines reproducing primes in series
(along with the continuation of the prime series) and
satisfying in addition the conditions b) and c).

A spline formed by polynomials with integer coefficients will
be called {\it the arithmetic spline}.

This article discusses two candidates for the arithmetic splines,
cubic one and parabolic.

These splines do not approximate a prime series.
Primes  are implanted in the structure of the splines,
which ensures that they are exactly reproducible. Such splines
lead to explicit soluble systems of linear equations whose
coefficients represent arithmetic functions themselves.

The inverse function  $p^{-1}(x)$ constitutes a differentiable
analogue of the number--theoretic function
$\pi (x)=\sum\limits_{p\leq x}1$ which  is comparable with
the Riemann exact expression for $\pi(x)$ through the zeros
of the $\zeta$--function(~\cite{Riemann}, page 34).\\
\newpage

{\bf 2. Cubic spline}\\

Consider the spline
$$
S_{cub}(x)=\left\{
\begin{array}{l}
x+1,\quad 1\leq x\leq 1.5,\\[1.5mm]
c_i (x),\quad i-0.5\leq x\leq i+0.5,\quad i=2, 3,\ldots,\\[1.5mm]
c_i(x)=c_{i+1}(x),\quad x=i+0.5,\quad i=1,2,\ldots,\\[2.5mm]
\frac{\di dc_i(x)}{\di dx}=\frac{\di dc_{i+1}(x)}{\di dx},
\quad x=i+0.5,\quad i=1,2,\ldots,
\end{array}
\right.
$$
with
$$
c_i(x)=2\left( a_i(x-i-0.5)^2 +b_i(x-i-0.5)+\frac{p(i)+p(i+1)}{2}
\right)(x-i)-
$$
$$
-2p(i)(x-i-0.5).
$$
Exact reproducibility of the primes follows from the
identity
\begin{equation}
\label{rone}
c_i(i)\equiv p(i),\quad i=1,2,\ldots\,.
\end{equation}
At the points of sewing the spline should also obey the identity
\begin{equation}
\label{rtwo}
c_i(i+0.5)\equiv \frac{1}{2}(p(i)+p(i+1)),\quad i=1,2\ldots ,
\end{equation}
which brings into the spline  additional information
of the prime series behaviour.

There exists  {\it an unique cubic spline} of the kind
$S_{cub}(x)$ with the coefficients
$$
\left.
\begin{array}{l}
{\di a_i=\frac{1}{2}(p(i+1)-p(i-1))-1,}\\
\phantom{fff}\\
b_i=p(i+1)-p(i)-1,
\end{array}
\right\} \quad i=2,3\ldots\, .
$$

This spline can be considered only as  almost--arithmetic.
The coefficients $\gamma_i$ and $\delta_i$ in
$c_i(x)=\alpha_i x^3+\beta_i x^2+\gamma_i x+\delta_i$
appear for some $i$ in the form $q+1/2,\, q\in N$.

\begin{figure}
\centering
\includegraphics[angle=-90, width=14cm, keepaspectratio]{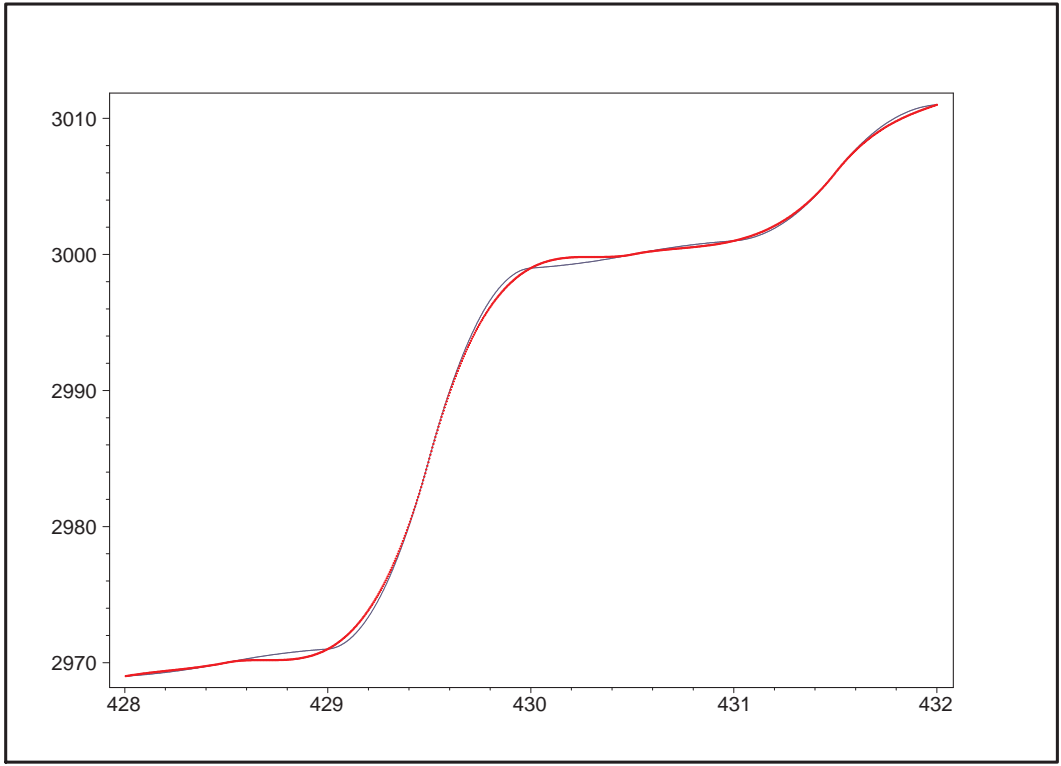}
\caption{ Functions $p(x)$: thick (red) line
corresponds to $S_{cub}(x)$;
 thin (blue) line  to $S_{quad}(x)$.}
\label{fig1}
\end{figure}

The negative value of the  discriminant
$$
d_i=4(p(i))^2-4(p(i-1)+p(i+1))p(i)+\frac{1}{4}(p(i-1))^2+\frac{1}{4}(p(i+1))^2+
$$
$$
+\frac{7}{2}p(i-1)p(i+1)+3
$$
gives positivity of the derivative
$$
\frac{\di dc_i(x)}{\di dx}=2(2a_i(x-i-0.5)+b_i)(x-i)+2a_i(x-i-0.5)^2
+2b_i(x-i-0.5)-
$$
$$
-p(i)+p(i+1).
$$

The inequality  $d_i<0$ leads to the following condition for
prime triplets\\ $p(i-1),\, p(i),\, p(i+1)$:
\begin{equation}
\label{rthree}
\frac{1}{2}(p(i-1)+p(i+1))-\frac{1}{4}\sqrt{t_i}<p(i)<\frac{1}{2}
(p(i-1)+p(i+1))+\frac{1}{4}\sqrt{t_i},
\end{equation}
where
$$
t_i=3((p(i+1)-p(i-1))^2-4)>0,\quad i=1,2,\ldots\, .
$$

Condition (\ref{rthree}) can be violated in the cases
\begin{equation}
\label{rfour}
\left.
\begin{array}{l}
p(i)-p(i-1)=\Delta_1,\quad p(i+1)-p(i)\geq \Delta_2;\,\\
\phantom{ss}\\
p(i)-p(i-1)\geq \Delta_2,\quad p(i+1)-p(i)= \Delta_1,
\end{array}
\right\}
\end{equation}
$$
\mbox{with \quad }\Delta_1=2,\quad \Delta_2=28.
$$

\begin{figure}
\centering
\includegraphics[angle=-90, width=14cm, keepaspectratio]{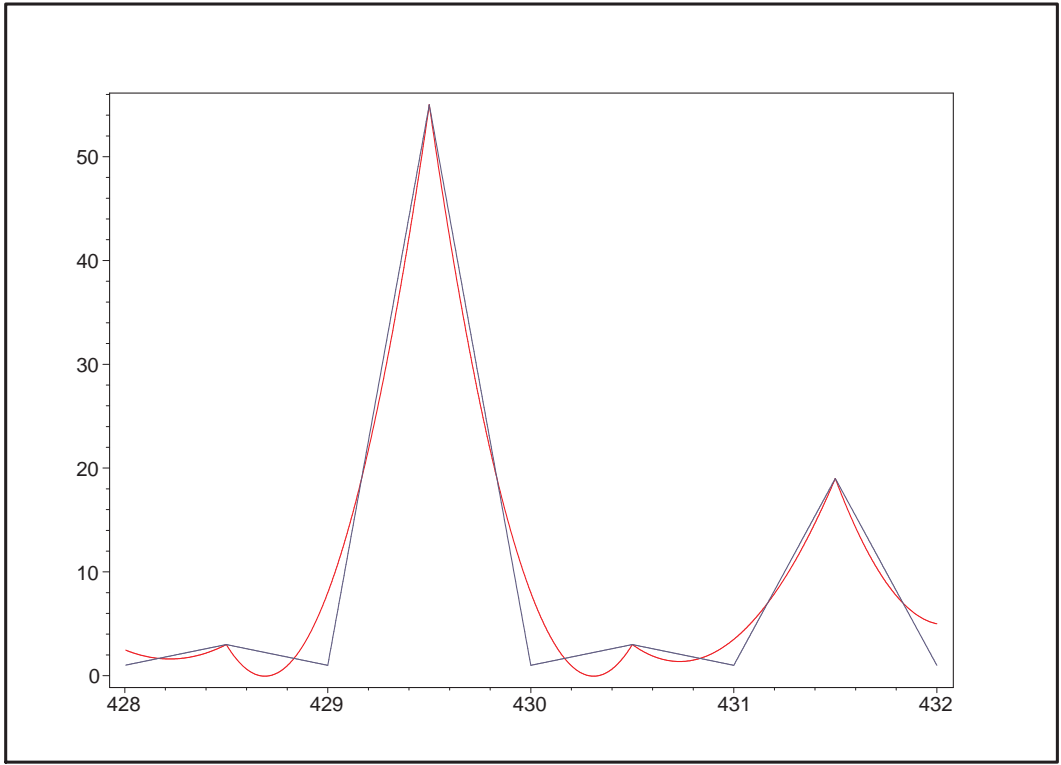}
\caption{Derivatives
$S'_{cub}(x)$(quadrics) and $S'_{quad}(x)$(segments)}
\label{fig2}
\end{figure}

Among the  first 1000 primes only 5 triplets
violate the rule (\ref{rthree}):
\begin{equation}
\label{rfive}
\begin{array}{l}
(2969,\, 2971,\, 2999),\\
(2971,\, 2999,\, 3001),\\
(3271,\, 3299,\, 3301),\\
(6917,\, 6947,\, 6949),\\
(7757,\, 7759,\, 7789).
\end{array}
\end{equation}

Except  these triplets  (including twin pairs),
condition (\ref{rthree}) is violated by
triplets of the  kind (\ref{rfour}) at the following values for
$\Delta_1$ and $\Delta_2$ :
$$
\begin{array}{l}
\Delta_1=4,\quad \Delta_2=56,\\
\Delta_1=6,\quad \Delta_2=84,\\
\Delta_1=8,\quad \Delta_2=114,\\ \mbox{ and all that.}
\end{array}
$$

Despite the cases where condition  (\ref{rthree}) is violated,
spline $S_{cub}(x)$ is convenient for creating subroutines
$p(x),\, \frac{\di dp(x)}{\di dx}$ and $p^{-1}(x)$, since the
inverse function
$p^{-1}(x)$ exists in the neighborhood of each prime number.\\

{\bf 3. Parabolic spline}\\

Given the following pairs of parabolas
$$
q_{i}(x)=\left\{
\begin{array}{l}
q_i^{l}(x),\quad i-0.5\leq x\leq i,\quad i=2, 3,\ldots,\\[1mm]
q_i^{r}(x),\quad i\leq x\leq i+0.5,\quad i=2, 3,\ldots,\\[3mm]
\left.
\begin{array}{l}
q_i^{l}(x)=q_{i}^{r}(x)\\[2mm]
\frac{\di dq_i^{l}(x)}{\di dx}=\frac{\di dq_{i}^{r}(x)}{\di dx}
\end{array}
\right\},\quad x=2,3,\ldots ; \mbox{ internal sewing,}
\end{array}
\right.
$$
with
$$
\begin{array}{l}
q_i^{l}(x)=-2a_{i-1}(x-i)^2 +x-i+p(i),\\[1mm]
q_i^{r}(x)=2a_i(x-i-0.5)^2+(2a_i+1)(x-i-0.5)+
\frac{\di p(i)+p(i+1)}{\di 2},\\[2mm]
a_i=p(i+1)-p(i)-1.\\[1mm]
\end{array}
$$

The parabolic spline
$$
S_{quad}(x)=\left\{
\begin{array}{l}
x+1,\quad 1\leq x\leq 1.5,\,\mbox{ initial polynomial,}\\[1mm]
q_i (x),\quad i-0.5\leq x\leq i+0.5,\quad i=2, 3,\ldots,\\[3mm]
\left.
\begin{array}{l}
q_i(x)=q_{i+1}(x)\\[2mm]
\frac{\di dq_i(x)}{\di dx}=\frac{\di dq_{i+1}(x)}{\di dx}
\end{array}
\right\},
\begin{array}{l}\quad x=i+0.5,\quad i=1,2,\ldots ;\\
\quad\mbox{external sewing}
\end{array}
\end{array}
\right.
$$
solves the problem better than the spline $S_{cub}$.
It has the following properties:

\noindent 1)  identities analogous to (\ref{rone}) and (\ref{rtwo})
are applicable

\begin{equation}
\label{pnew1}
q_i(i)\equiv p(i),
\end{equation}
\begin{equation}
\label{pnew2}
q_i(i+0.5)\equiv\frac{\di 1}{\di 2}(p(i)+p(i+1)),\quad i=1,2,
\ldots ;
\end{equation}

2) the derivatives
$$
\frac{\di dq_i^{l}(x)}{\di dx}=4a_{i-1}(i-x)+1,
$$
$$
\frac{\di dq_i^{r}(x)}{\di dx}=4a_i(x-i)+1
$$
take a minimal value  $+1$ at the points of internal sewing
(they are the points of interpolation to the spline) and maximal
values at the points of external sewing, where they coincide
with the derivative of the spline $S_{cub}$.

\begin{figure}
\centering
\includegraphics[angle=-90, width=14cm, keepaspectratio]{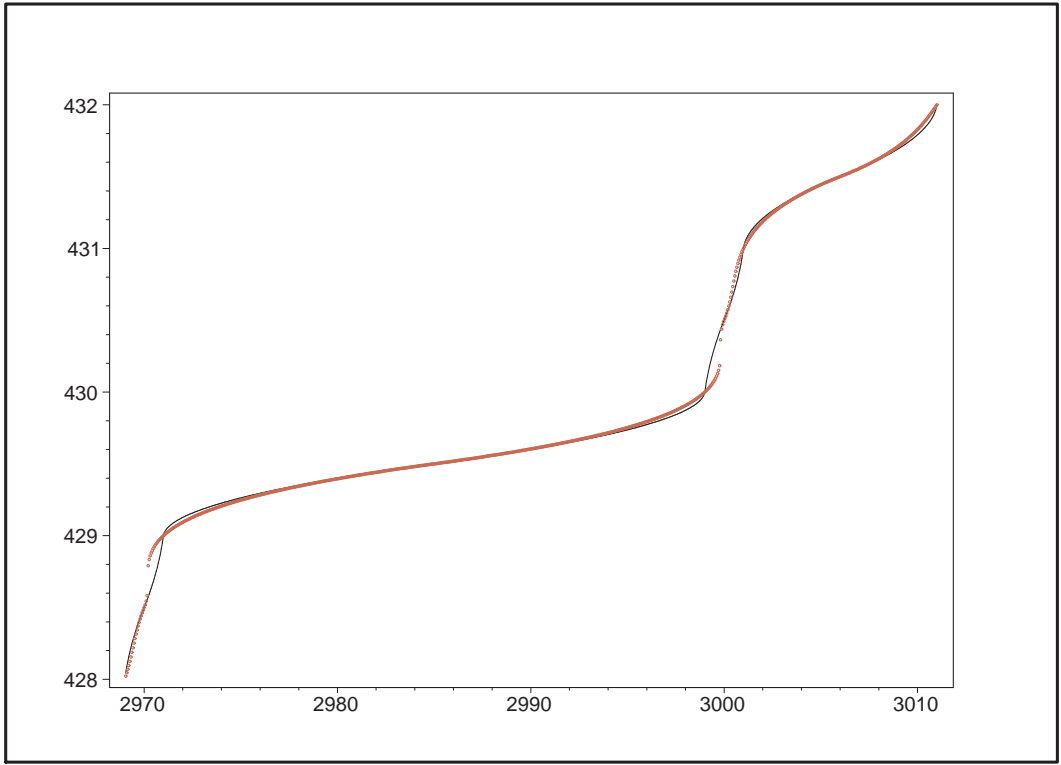}
\caption{Inverse functions
$p^{-1}(x)$: thick line corresponds to $S_{cub}$, thin
line to $S_{quad}$.}
\label{fig3}
\end{figure}

The positive values of the derivatives show that the spline
$S_{quad}(x)$ monotonically increases on the semi-axis
$[1,\,\infty)$.

There exists a function $S_{quad}^{-1}(x)$,
inverse to the function $S_{quad}(x)$, determined on the axis
$[2,\,\infty )$ and thus the spline $S_{quad}$ fulfills the conditions a),
b), c). Moreover, the function  $S_{quad}^{-1}(x)$ is differentiable
on $(2,\,\infty )$.

 The spline $S_{quad}$ is {\it arithmetic} because
 {\it prime number polynomials}
$$
\begin{array}{l}
q_i^{l}(x)=\alpha_i^{l} x^2+\beta_i^{l} x+\gamma_i^{l}\quad and \quad
q_i^{r}(x)=\alpha_i^{r} x^2+\beta_i^{r} x+\gamma_i^{r}\quad
\end{array}
$$
hold integer coefficients (see Table 1):

$$
\alpha_i^{l}=-2a_{i-1};
\quad\beta_i^{l}=4ia_{i-1}+1;
\quad\gamma_i^{l}=-2i^2a_{i-1}+p(i)-i;
$$
$$
\alpha_i^{r}=2a_i;
\quad\beta_i^{r}=-4ia_i+1;
\quad\gamma_i^{r}=2i^2a_i+p(i)-i.\\
$$
The joint satisfiability
of identities (\ref{pnew1}), (\ref{pnew2}) reinforces the
hypothesis that
{\it the spline}
 $S_{quad}$ {\it is an unique arithmetic spline satisfying
 the conditions} a), b) and c).\\[1cm]

\begin{center}
{\large {\bf Table 1}} $\left(d_i^{l}=(\beta_i^{l})^2-
4\alpha_i^{l}\gamma_i^{l},\quad d_i^{r}=
(\beta_i^{r})^2-
4\alpha_i^{r}\gamma_i^{r}\right)$\\
\end{center}

\noindent
$\underline{\phantom{bbb}i\hspace{0.6cm} p(i)\hspace{1.6cm}
 \alpha_i^{l}\hspace{0.7cm}
\beta_i^{l}\hspace{0.9cm}\gamma_i^{l}
\hspace{1.0cm}d_i^{l}
\hspace{1.8cm} \alpha_i^{r}
\hspace{0.9cm}\beta_i^{r}\hspace{0.9cm}\gamma_i^{r}
\hspace{1.0cm}d_i^{r}}$
\begin{verbatim}
  2    3         0     1      1      1          2     -7      9     -23
  3    5        -2    13    -16     41          2    -11     20     -39
  4    7        -2    17    -29     57          6    -47     99    -167
  5   11        -6    61   -144    265          2    -19     56     -87
  6   13        -2    25    -65    105          6    -71    223    -311
  7   17        -6    85   -284    409          2    -27    108    -135
  8   19        -2    33   -117    153          6    -95    395    -455
  9   23        -6   109   -472    553         10   -179    824    -919
 10   29       -10   201   -981   1161          2    -39    219    -231
 11   31        -2    45   -222    249         10   -219   1230   -1239
 12   37       -10   241  -1415   1481          6   -143    889    -887
 13   41        -6   157   -986    985          2    -51    366    -327
 14   43        -2    57   -363    345          6   -167   1205   -1031
 15   47        -6   181  -1318   1129         10   -299   2282   -1879
 16   53       -10   321  -2523   2121         10   -319   2597   -2119
 17   59       -10   341  -2848   2361          2    -67    620    -471
 18   61        -2    73   -605    489         10   -359   3283   -2439
 19   67       -10   381  -3562   2681          6   -227   2214   -1607
 20   71        -6   241  -2349   1705          2    -79    851    -567
 21   73        -2    85   -830    585         10   -419   4462   -2919
\end{verbatim}

It should be noted here that the  coefficient
$a_i$ in  $q_i^{l}(x)$ and $q_i^{r}(x)$
represents a basic arithmetic function --
{\it the number of composite numbers in the interval} $(p(i),
p(i+1))$.

Figures 1--3 present a comparison of the splines  $S_{cub}$ and
$S_{quad}$. The interval's $[428, 432]$ image (argument to the functions $p(x)$
and $dp(x)$) contains the first pair of triplets
 (\ref{rfive}) violative the positivity
 of the derivative $S'_{cub}(x)$.

Figure 3 shows intervals where no inverse function $p^{-1}(x)$ for the spline
$S_{cub}$ exists.

Figure 2 illustrates the properties of the derivatives $S'_{cub}$ and $S'_{quad}$:
reaching the minimal and maximal
values and the equalities of derivatives at the sewing
points.\\

{\bf 4. Inverse parabolic spline $S_{quad}^{-1}(x)$ and its derivative } \\

The pairs of functions

$$
t_{i}(x)=\left\{
\begin{array}{l}
t_i^{l}(x),\quad \frac{\di p(i-1)+p(i)}{\di 2}\leq x\leq p(i),
\quad i=2, 3,\ldots,\\[1mm]
t_i^{r}(x),\quad p(i)\leq x\leq \frac{\di p(i)+p(i+1)}{\di 2},
\quad i=2, 3,\ldots,\\[3mm]
\left.
\begin{array}{l}
t_i^{l}(x)=t_{i}^{r}(x)\\[2mm]
\frac{\di dt_i^{l}(x)}{\di dx}=\frac{\di dt_{i}^{r}(x)}{\di dx}
\end{array}
\right\},\quad x=p(2), p(3),\ldots ,
\end{array}
\right.
$$
where

$$
\begin{array}{l}
t_i^{l}(x)=
i+\frac{\di 1-\sqrt{b_i^{l}}}{\di 4a_{i-1}},
(t_2^{l}=x-1),\quad
t_i^{r}(x)=i+
\frac{\di \sqrt{b_i^{r}}-1}{\di 4a_i},\\[4mm]
b_i^{l}=8a_{i-1}(p(i)-x)+1,\quad b_i^{r}=8a_{i}(x-p(i))+1,\\[3mm]
a_i=p(i+1)-p(i)-1,
\end{array}
$$

\vspace{0.5cm}

determine the {\it inverse spline}

$$
S_{quad}^{-1}(x)=\left\{
\begin{array}{l}
x-1,\quad 2\leq x\leq 2.5,\\[1mm]
t_i(x),\quad \frac{\di p(i-1)+p(i)}{\di 2}\leq x\leq
\frac{\di p(i)+p(i+1)}{\di 2},\, i=2, 3,\ldots,\\[3mm]
\left.
\begin{array}{l}
t_i(x)=t_{i+1}(x)\\[2mm]
\frac{\di dt_i(x)}{\di dx}=\frac{\di dt_{i+1}(x)}{\di dx}
\end{array}
\right\},  x=\frac{\di p(i)+p(i+1)}{\di 2},
\, i=1, 2,\ldots\, .\\[3mm]
\end{array}
\right.
$$

The derivative of the inverse spline is as follows

$$
\frac{\di dS_{quad}^{-1}(x)}{\di dx}=\left\{
\begin{array}{l}
1,\quad 2\leq x\leq 2.5,\\[2mm]
\left(b_i^{l}\right)^{-1/2},\quad
\frac{\di p(i-1)+p(i)}{\di 2}\leq x\leq p(i),\,
i=2, 3,\ldots,\\[3mm]
\left(b_i^{r}\right)^{-1/2},\quad
p(i)\leq x\leq \frac{\di p(i)+p(i+1)}{\di 2},\,
i=2, 3,\ldots\, .
\end{array}
\right.
$$

\begin{figure}
\centering
\includegraphics[angle=-90, width=14cm, keepaspectratio]{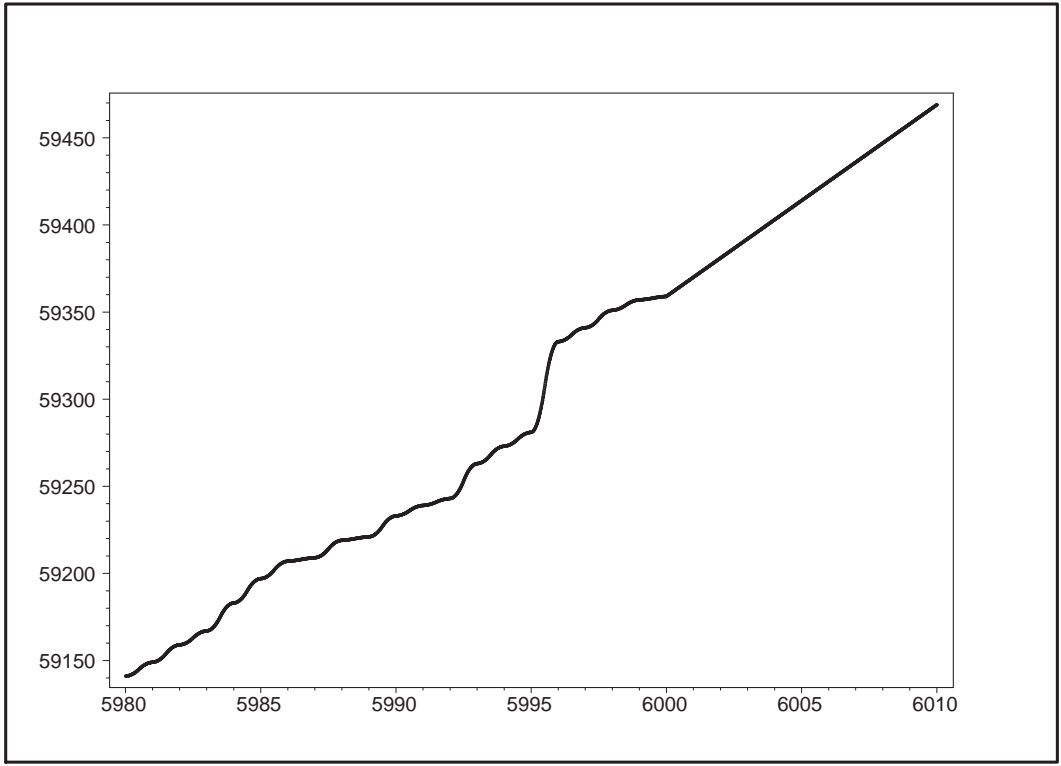}
\caption{ Sewing the spline $S_{quad}(x)$ with asymptote
$\widetilde{p}(x)$ at  x=6000.}
\label{fig1}
\end{figure}

 {\bf 5. About subroutines   $p(x),\,
 \frac{\di d p(x)}{\di dx}$, $p^{-1}(x)$} and $\frac{\di dp^{-1}(x)}{\di dx}.$\\

Both splines  $S_{cub}$ and $S_{quad}$ were employed to create
Fortran functions $p(x)$,\linebreak $dp(x):=\frac{\di dp(x)}
{\di dx}$,
$p_{-}(x)$, $p_{-}newt(x)$ ($p_{-}(x)$ and
$p_{-}newt(x)$ denote $ p^{-1}(x)$) and\linebreak $dp_{-}(x)
:=\frac{\di dp^{-1}(x)}{\di dx}$.

For
convenience in applications the functions $p(x)$ and $dp(x)$ have
been extended to  $+\infty$ by the  asymptote (\cite{four},
page 140)\\
\begin{equation}
\label{rp}
\widetilde{p}(x)=x\left(\ln x +\ln\ln x+\frac{\ln\ln x-2}{\ln x}-
\frac{{(\ln\ln x)}^2/2-3\ln\ln x+5.5}{{(\ln x)}^2}-1\right).
\end{equation}

To obtain primes, two alternative programs have been
employed -- {\it subroutine eratosthenes(n)} and
 {\it subroutine primes(n)}.
The former accomplishes this by generating primes up to a
given value $n$, whereas the latter achieves the result by
reading $n$ primes from  given 6 column file named
{\it primes}.  In these programs, the spline and its derivative
are automatically
sewed with asymptote (\ref{rp}) and its derivative.

For purposes of building the inverse function $p^{-1}(x)$ in
the programs $p_{-}newt(x)$ and $p_{-}(x)$, two different
approaches have been used.

The program $p_{-}newt(x)$ is based on an {\it autoregularized
variant of the Newton method} (\cite{alex1}, page 43).\\

\begin{equation}
\label{cc}
\left.
\begin{array}{l}
y_{0},\quad \varepsilon_0>0,\quad y_{k+1}=y_k-\frac{\di p(y_k)-x}{\di dp(y_k)+\varepsilon_k},\quad
k=0,1,2,\ldots,\\
\phantom{dcd}\\
\varepsilon_k=\frac{\di 1}{\di 2}\left(\sqrt{(dp(y_k))^2+4N|p(y_k)-x|}-
dp(y_k) \right),\\
\phantom{dcd}\\
N= (\varepsilon_0^2+\varepsilon_0 dp(y_0))/|p(y_0)-x|,
\end{array}
\right\}
\end{equation}
where several combinations of the initial value $y_0=li(x)$ and
initial regularizator  $\varepsilon_0$ (these combinations can be
seen at the beginning of the $p_{-}newt(x)$ program's   body)
ensure a construction of the inverse function $p^{-1}(x)$ on the
interval $[2, 10^8].$

The way of setting the initial approximation $y_0$ within the
values
$$
\frac{\di x}{\di \ln x},\quad li(x)=\int\limits_{2}^{x}\frac{\di dt}{\di \ln t},
$$
\begin{equation}
\label{riem1}
R(x)=
\sum\limits_{n=1}^{\infty}\frac{\mu (n)}{n}f(x^{1/n})\quad
\mbox{(see \cite{Riemann},\, page 35),}
\end{equation}
where  $f(x)=li(x)$, $\mu(n)$ is $0$ if $n$ is divisible by  a prime square,
$1$ if $n$ is a product of an even number
of distinct primes, and $-1$ if $n$ is a product  of an odd
number of distinct primes
was checked: the conclusion is that
the method (\ref{cc}) and $y_0=li(x)$ are an acceptable combination.

Figure 4 shows the sewing of the function $p(x)$ with the asymptote
(\ref{rp}).

The programs $p_{-}(x)$ and $dp_{-}(x)$  are based on the application
of the inverse spline $S^{-1}_{quad}$ and its derivative. In these
programs, an algorithm has been applied by virtue of which the
needed pairs $t_{i}(x)$ and $b_{i}(x)$ are established by
approximating $ix = \lfloor li(x)\rfloor$ (see the beginning of the programmes
$p_{-}(x)$ and $dp_{-}(x)$). It is worth mentioning that this
algorithm works in the case where $\pi(x) < li(x)$, as well as in
the case where $\pi(x) > li(x)$, i.e., the algorithm does not
depend on the knowledge the minimal value of $x$ for which
the difference $li(x)-\pi(x)$ changes the sign.

The programs $p(x)$, $dp(x)$, $p_{-}newt(x)$, $p_{-}(x)$ and
$dp_{-}(x)$ have been
realized in the Fortran90.
These programs form the package named
$pp_{-}f90[1, \infty)$, available in Appendix 1. They employ the
{\it natural} extension of the function to $-\infty$, based on
the values  of the initial polynomial $x+1$.

The package $pp_{-}f90[1, \infty)$ is immediately applicable to
{\it Compaque}-- and {\it MS-- Fortran} and is
facile transportable to other Fortran versions.

Fortran functions $p(x)$ and  $p^{-1}(x)$ are as easily
applicable as the intrinsic functions  $sin(x)$ and $exp(x)$.

The end of Appendix 1 contains a program called $test_{-}pp_{-}$
 which has
been used to compute all the tables supporting the graphics in
this article and which serves as an illustration
for application of the functions  $p(x)$ and  $p^{-1}(x)$.\\

{\bf 6. Possible applications of the functions $p(x)$ and  $p^{-1}(x)$.}\\

{\bf 6.1.} Functions  $p(x)$ and $p^{-1}(x)$ can be used to introduce
new functions $\sin p(x)$, $\cos p(x)$, $\tan p(x)$,
$e^{-p(x)}$ and $\ln p(x)$,
applicable when the specific character of nonasymptotic
prime number distribution  is necessary to be accounted for
(see, e.g., \cite{eight}).
In particular, the functions $\sin p(x)$ and $\cos p(x)$
can be used for creation of a {\it prime number harmonic analysis}.\\

{\bf 6.2.} Diophantine equations can be solved within the
following approximate method: to the Diophantine equation
\begin{equation}
\label{diof1}
f_1(x_1, ..., x_n)=0
\end{equation}
one adds a new equation, {\it from reals-to-integers equation}
\begin{equation}
\label{diof2}
f_2^{h}(x_1, ..., x_n)
: =\sin^2 (\pi x_1)+\sin^2 (\pi x_2)+\cdots +\sin^2 (\pi x_n)=0,
\end{equation}
or one adds {\it from reals-to-primes equation}
\begin{equation}
\label{diof3}
f_2(x_1, ..., x_n)
: =\sin^2 (\pi p^{-1}(x_1))+\sin^2 (\pi p^{-1}(x_2))+\cdots +
\sin^2 (\pi p^{-1}(x_n))=0.
\end{equation}

\begin{figure}
\centering
\includegraphics[width=18cm, height=16cm, angle=-90]{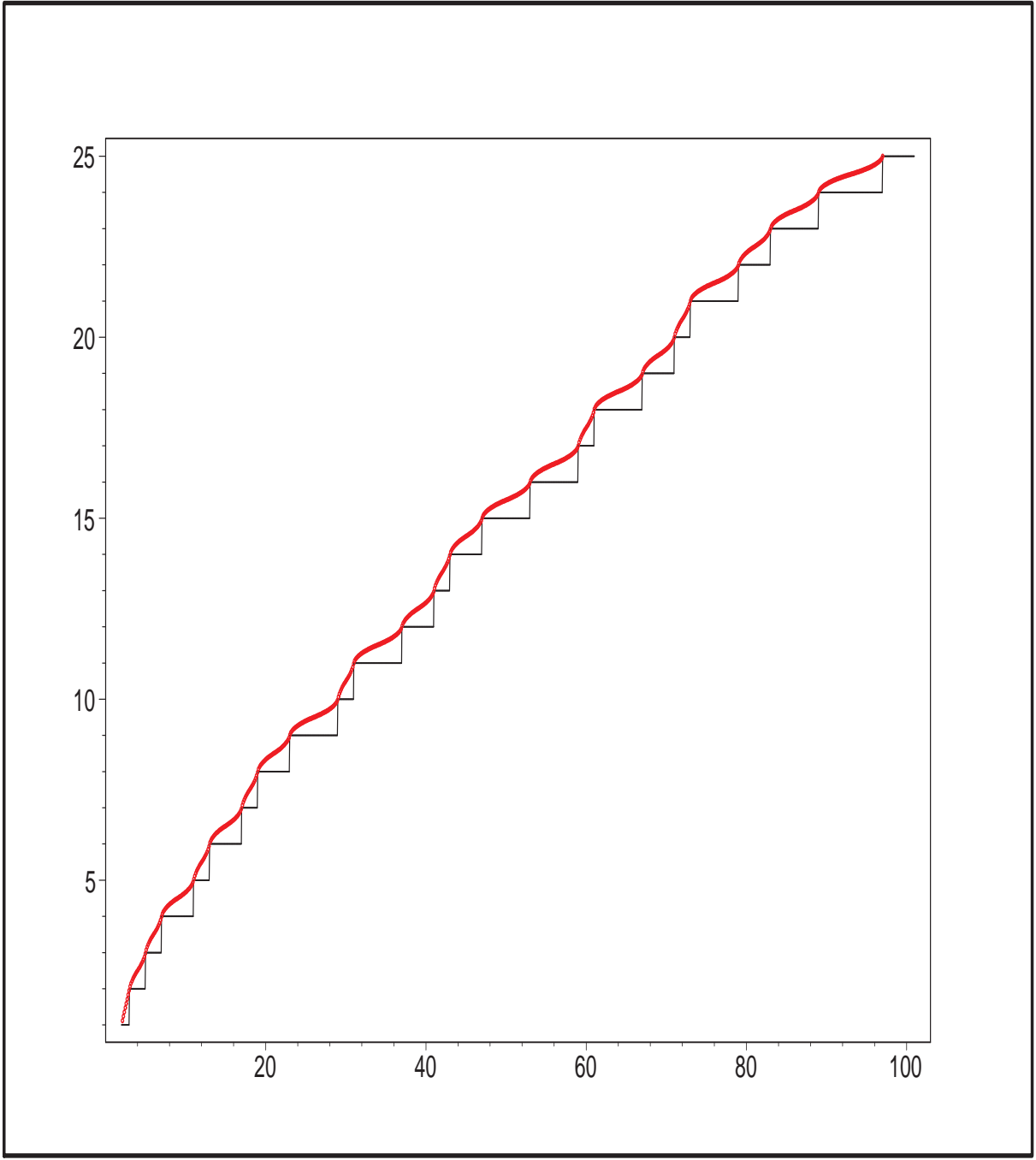}
\caption{
The differentiable inverse function $p^{-1}(x)$ (thick line)
versus the Riemann--Von Mangoldt step--function $\pi_{R}(x)$.}
\end{figure}

By solving either system  (\ref{diof1}), (\ref{diof2}) or
(\ref{diof1}), (\ref{diof3}) one can find solutions to
(\ref{diof1}) as real  approximations to the natural or
the prime numbers.

The above systems can be solved by methods
working in the case of degeneration
of the derivative at the solution (see, e.g., \cite{nine}).

Here two examples of solving such systems by means of autoregularized
iterative processes ($rgn$) \cite{five} are presented;
$rgn$-processes are combined
with both the $svd$-method \cite{gol} and the adaptive scaling \cite{mor}.

In this case the program {\it afxy} \cite{alex2} is used to
find all solutions of the systems  (\ref{diof1}), (\ref{diof2}) and
(\ref{diof1}), (\ref{diof3}) in a given definition domain.

Let the problem (\ref{diof1}), (\ref{diof3}) be written as
\begin{equation}
\label{diof4}
Fx=0,
\end{equation}
where
$$
Fx=f'^{T}(x)(fx-\overline{y}),\, fx=[f_1(x),\, f_2(x)]^{T},\,
f:D_f\subset R^n \rightarrow R^m,\, x\in R^n,\,
\overline{y}\in R^m,
$$
and $D_f$ is an
open convex domain in $R^n$.

The linear problem at the $k$th iteration of $rgn$-process is of the kind
\begin{equation}
\label{diof5}
(f'^{T}(x^{k})f'(x^{k})+\varepsilon_k I)(x^{k+1}-x^{k})
=-Fx^{k},
\end{equation}
where
$$
\varepsilon_k=\frac{1}{2}\left(\sqrt{\tau_k^2+4N\rho_k}-\tau_k
\right),
$$
$$
\tau_k=\|f'^{T}(x^{k})f'(x^{k})\|_{\infty},\,
\rho_k=\|Fx^{k}\|_{\infty},\,
 N=(\varepsilon_0+\varepsilon_0\tau_0)/\rho_0.
$$
For simplicity in equality (\ref{diof5}) and in the expressions
for $\tau_k, \rho_k$ and $\varepsilon_k I$, the scaling operators are
neglected.
Just for the problem  (\ref{diof5}) the $svd$-method is in use.

For the purpose of finding all solutions of equation   (\ref{diof4})
in $D_f$, the program {\it afxy} realizes an algorithm in which the vector
$Fx_k$ is multiple factored by local {\it extractors} of the kind
$$
e_r(x, \overline{x}^{(r)})=\left(1-e^{-\| x-\overline{x}^{(r)}\|_{_2}}
\right)^{-1},
$$
where $\overline{x}^{(r)}\in D_f$ is the $r$th solution of equation
(\ref{diof4}). The transformed problem
\begin{equation}
\label{diof6}
F_{r^*}x: = \left(\prod\limits_{r=1}^{r^*}e_r(x,
\overline{x}^{(r)})\right) Fx=0,\quad r^*
\geq 1
\end{equation}
is multiple solvable by the program {\it afxy}.

For each new problem (\ref{diof6}) {\it afxy} realizes different $rgn$ iterative
processes: with different guesses $x_0$ (randomly formed by the initially
 given $x_0$ and by peculiarities of the domain $D_f$) and with different
initial regularizators $\varepsilon_0$ (from the preset table of
regularizators).\\

{\bf Example 1.}
Solution of the equation $x_1^2+x_2^2=x_3^2+1$ over primes.

Consider problem  (\ref{diof4}) with
$$
fx=\left\{
\begin{array}{l}
f_1(x): =x_1^2+x_2^2-x_3^2-1=0,\\
\phantom{r}\\
f_2(x): =\sin^2 (\pi p^{-1}(x_1))+\sin^2 (\pi p^{-1}(x_2))+\cdots +
\sin^2 (\pi p^{-1}(x_n))=0,
\end{array}
\right.
$$
$$
n=3, m=2, D_f=[2, 100]\times[2, 100]\times[2, 100], x=(x_1, x_2, x_3)^T, \,
\overline{y}=[0, 0]^T.
$$
{\it Run-time section} of Application 2 contains a
subroutine $fxy(m, n, np, neq, f, x, pp, df, yr)$, where the equation $fx=0$
and derivatives  $\frac{\di df_j(x)}{\di dx_i}\,  (i=1,2,3,\, j=1,2)$
are coded. Here, in {\it the pre-exe section} of the program {\it afxy},
the main controls
$m, n, np, f, lsmh, mqh, nsolh,$ $x_1^{(0)},x_2^{(0)}$ and $D_f$
are given as well.

If in {\it the pre-exe section} some
main control is not prescribed, then {\it afxy}
switches to an interactive mode and demands from display an adjustment to
the value of this control.

The solution $(x_1, x_2, x_3)$ of the equation $f_1(x)=0$
(\cite{serp}, page 35)
is thought  to be a {\it quasi-Pythagorean prime triplet } (the equation
$x_1^2+x_2^2=x_3^2$ has no solutions on primes). There exist two
series of natural numbers satisfying equation $f_1(x)=0:$
$$
\left(2n+1,\, n^2+n-1,\, n^2+n+1\right),\,
\mbox{\cite{serp},\, page 35},
$$
$$
\left(2n(4n+1),\, 16n^3-1,\, 16n^3+2n\right),\,
\mbox{\cite{serp},\, page 36}.
$$

The present state of the program {\it afxy} can only produce
up to 20 (i.e. $r^*\leq 20$)
extractions and can find only 20 {\it
quasi-Pythagorean prime triplets} in $D_f$(see FOUND SOLUTIONS in Application 2).\\

{\bf Example 2.} Solution of the equation $x_1^2+x_2^2=x_3^2+1$ over twins.\\

Consider problem  (\ref{diof4}) with
$$
fx=\left\{
\begin{array}{l}
f_1(x): =x_1^2+x_2^2-x_3^2-1=0,\\
\phantom{r}\\
f_2(x): =\sin^2 (\pi p^{-1}(x_1))+\sin^2 (\pi p^{-1}(x_2))+\cdots +
\sin^2 (\pi p^{-1}(x_n))=0,\\
\phantom{r}\\
f_3(x): =x_3-x_1-2=0,\\
\end{array}
\right.
$$
$$
n=3, m=3, D_f=[2, 100]\times[2, 100]\times[2, 100], \mbox{ and }\,
\overline{y}=[0, 0, 0]^T.
$$

The needed subroutine {\it fxy} for this example is presented in Application 3.

Application 3 shows that the program {\it afxy} finds all 5
{\it quasi-Pythagorean prime triplets}
containing twin pairs in the domain $D_f$.

{\bf 6.3.} In Figure 5, the inverse function  $p^{-1}(x)$ is
compared with the step--function $\pi (x)$
and with the Riemann--Von
Mangoldt continuous step--function $\pi_R (x)$ which results from
(\ref{riem1}) by means of the substitution
$$
f(x)=li(x)-\sum\limits_{n=1}^{\infty}li(e^{\di \,\rho_n\ln x})+
\int\limits_{x}^{\infty}\frac{\di da}{\di (a^2-1)a\ln a}-\ln 2,
$$
where $\{\rho_n\}$ are the complex zeros of the equation

\begin{equation}
\label{riem2}
\zeta(s)=0
\end{equation}
in the form $\rho_n=\frac{\di 1}{\di 2}\pm it_n,\, n=1,2,\ldots\, $,
and $\zeta(s)$ is Riemann's $\zeta$--function \cite{Riemann}.

In our case, counting function  $\pi (x)$ is expressed by the
inverse function $p^{-1}(x)$ as the formula

\begin{equation}
\pi(x)=\lfloor p^{-1}(x)\rfloor.
\end{equation}

Comparing the two ways, the function  $p^{-1}(x)$ proves more
convenient for application than the Riemann--Von Mangoldt function
$\pi_R (x)$; moreover,
$p(x)$ and  $p^{-1}(x)$ are applicable just now, not waiting for
the final representation of the function $\pi_R (x)$ through the zeros $\rho_n$.

{\bf 6.4.}For purposes of investigation the nonasymptotic
behaviour of primes, instead of the functions
$p(x)$ and $p^{-1}(x)$, one can resort the functions:

$$
A(x)=p(x)-\widetilde{p}(x)-(p(x_0)-\widetilde{p}(x_0)),\quad
x\in [x_0, x_0+\varepsilon(x_0))
$$
{\it --a local variance of the function} $p(x)$, where
$\varepsilon (x_0)$ is relatively small in comparison with $x_0$,
and

$$
B(x)=p^{-1}(x)-R(x)-(p^{-1}(x_0)-R(x_0)),\quad
x\in [x_0, x_0+\varepsilon(x_0))
$$
{\it --a local variance of the function $p^{-1}(x)$}.\\
Here $\widetilde{p}(x)$ is the asymptote (\ref{rp}), and $R(x)$
is the Riemann simplified formula (\ref{riem1}) for  $\pi (x)$.
Figures 7-9 serve as examples to the behaviuor of $A(x)$ and $B(x)$.\\

\begin{figure}
\centering
\includegraphics[angle=-90, width=13cm, keepaspectratio]{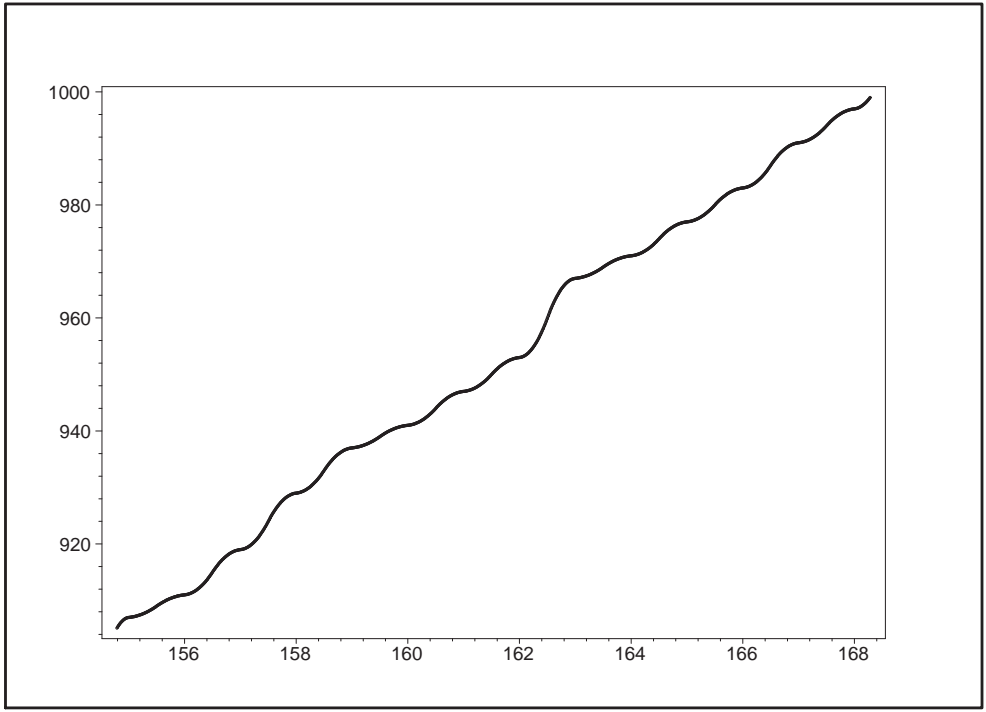}
\caption{{\it Function  $p(x)$.}
The interval [154.78, 168.2]  is chosen so as the definition
interval [900, 1000] of the function $p^{-1}(x)$ coincides
with the definition interval of the function
{\it number variance of the zeros $\rho_n$} from Figure 4
\cite{ten},\, p. 406.}
\label{fig6}
\end{figure}

\begin{figure}
\centering
\includegraphics[angle=-90, width=13cm, keepaspectratio]{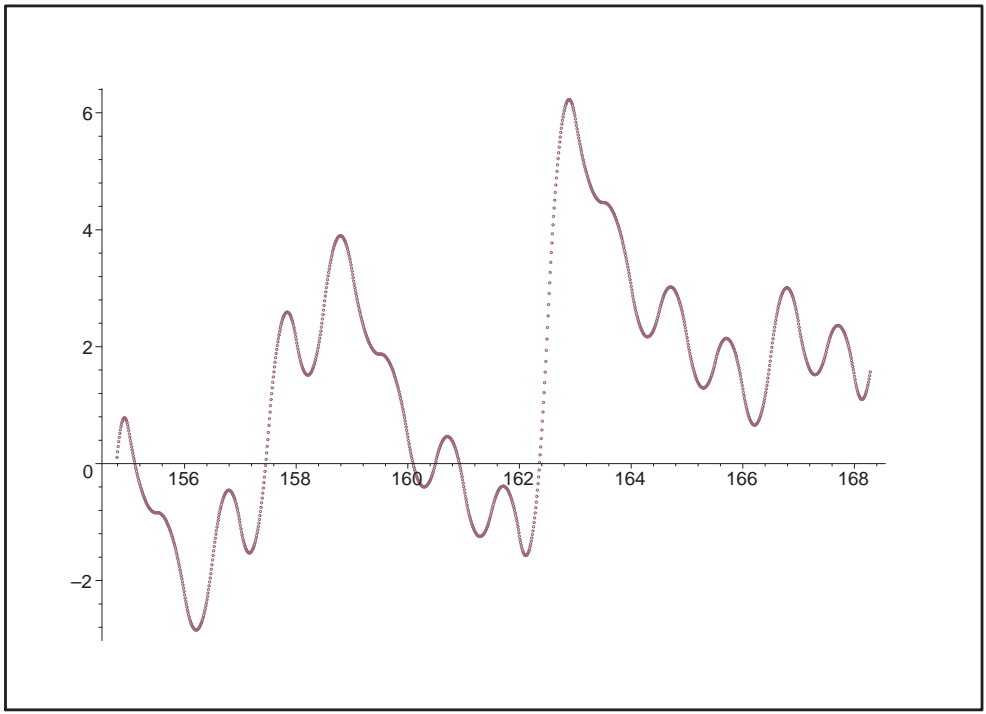}
\caption{{\it Function $A(x)$;}
note {\it a qualitative intimacy} with the function
{\it  number variance of the zeros $\rho_n$} from Figure 4
\cite{ten},\, p. 406, as well as the same number of peaks equal to 14
in the given interval; the function $A(x)$ covers 14 successive primes:
907, 911, 919, 929, 937, 941, 947,953, 967, 971, 977, 983, 991 and 997.}
\label{fig7}
\end{figure}

\begin{figure}
\centering
\includegraphics[angle=-90, width=13cm, keepaspectratio]{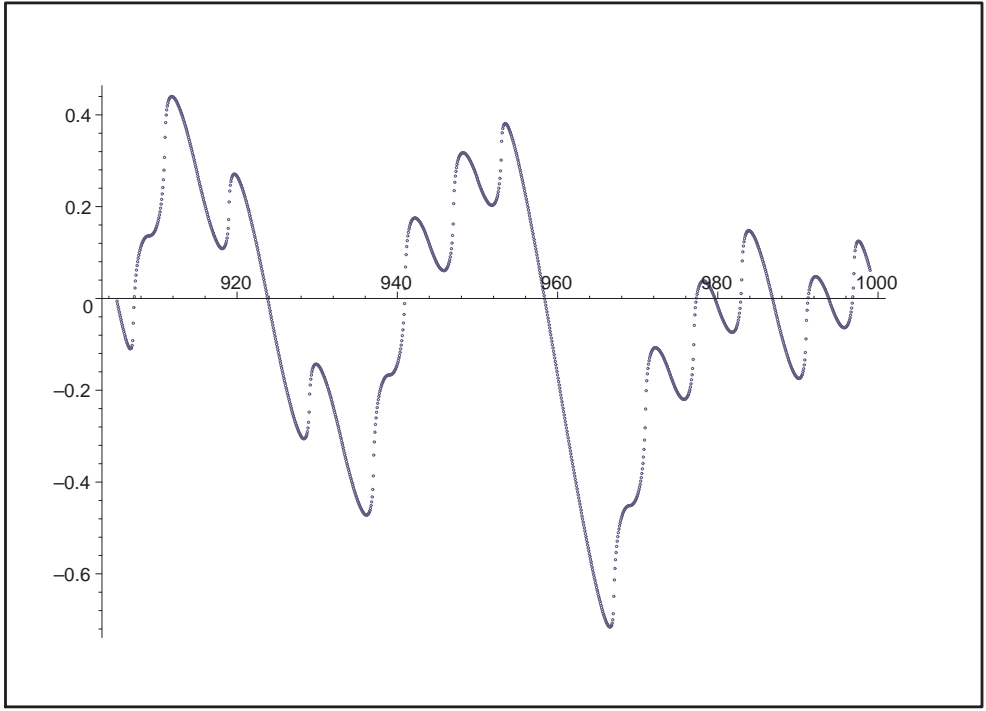}
\caption{{\it Function $B(x)$.}
{\it Inverse function} to the function $A(x)$ from Figure 7.}
\label{fig8}
\end{figure}

\begin{figure}
\centering
\includegraphics[angle=-90, width=13cm, keepaspectratio]{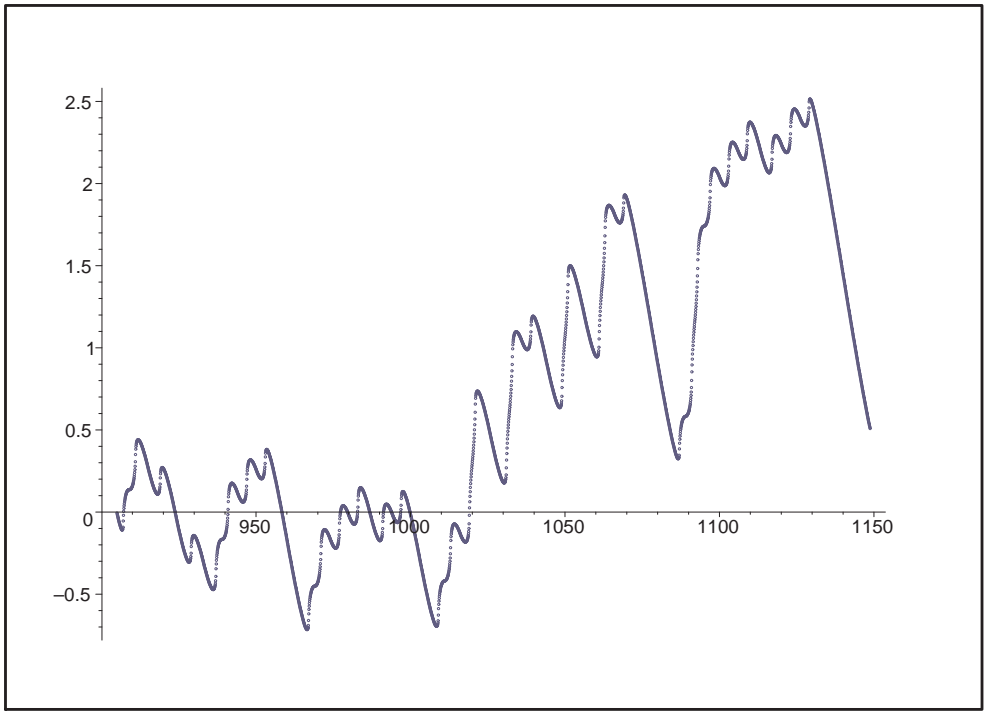}
\caption{{\it Function $B(x)$}
on the enlarged interval [900,1150]. How many peaks does the function
{\it number variance of the zeros $\rho_n$} have on the
interval [1000, 1150]? Do we really have no 17 peaks?}
\label{fig9}
\end{figure}

\newpage

\large
{\bf Application 1.} Programm package $pp_{-}f90[1, \infty)$\\
\small
\begin{verbatim}
       function p(x)   ! 1 <= x < infinity 
       implicit real*8(a-h,o-z)
       common/wn/wn
        if(x  <=  2.d0) then
           p = x+1.d0; return
        endif
        if(x > 2.d0 .and. x <= wn) then
           ix=floor(x+0.5d0); xx=dfloat(ix); p= sq(ix,x,xx); return
        endif
        if(x > wn) then
           p= r(x); return
        endif
        contains
        function sq(ix,x,xx)  ! left-right quadric pair
        common/protiarithmi/q(10000000)
        if(x <= xx) then
           sq=-2.d0*(q(ix)-q(ix-1)-1.d0)*(x-xx)**2+x-xx+q(ix)
        else
           sq=2.d0*(q(ix+1)-q(ix)-1.d0)*(x-xx-0.5d0)**2+            &
        (2.d0*(q(ix+1)-q(ix))-1.d0)*(x-xx-0.5d0)+(q(ix)+q(ix+1))/2.d0
        endif
        end function sq
        end function p

       function dp(x)  ! 1 <= x < infinity ; derivative of p(x) 
       implicit real*8(a-h,o-z)
       common/wn/wn
        if(x  <=  2.d0) then
           dp = 1.d0; return
        endif
        if(x > 2.d0 .and. x <= wn) then
           ix=floor(x+0.5d0); xx=dfloat(ix); dp= dsq(ix,x,xx); return
        endif
        if(x > wn) then
           dp=dr(x); return
        endif
        contains
        function dsq(ix,x,xx) ! derivative of left-right quadrics
        common/protiarithmi/q(10000000)
        if(x <= xx) then
           dsq=-4.d0*(q(ix)-q(ix-1)-1.d0)*(x-xx)+1.d0
        else
           dsq=4.d0*(q(ix+1)-q(ix)-1.d0)*(x-xx-0.5d0)+               &
               2.d0*(q(ix+1)-q(ix))-1.d0
        endif
        end function dsq
        end function dp
!
!
!
!
!
       function p_newt(u)      ! 2 <= x < infinity 
       implicit real*8(a-h,o-z)
       parameter(xy0=1.d3,xeps1=58.d3,eps1=3.d0,xeps2=17.d5,   &
       eps2=300,eps3=700.d0,ytol=1.d-11,dytol=1.d-11,ktol=1000)
       common/rcorr/rc,drc/kk/k/rw/r
       drcw=drc; drc=0.d0
       x=u; xw=x
       if(x >= 3.d0) then
          if(x < xy0) then
             y=x/dlog(x)
          else
             y=dlogintegral(x)
          endif
       else
          p_newt=x-1.d0
          return
       endif
       if( x < xeps1)  eps0=eps1           ! initial value of
       if(xeps1 <= x < xeps2)  eps0=eps2   ! initial value of
       if(xeps2 <= x)  eps0=eps3           ! autoregularizator
20     dt0=dp(y); r0=dabs(p(y)-x)          ! constant for
       en=(eps0**2+eps0*dabs(dt0))/r0      ! autoregularizator formula
       r=0.d0; k=0; small=1.d200
1      k=k+1
       yy=y; rr=r; t=p(yy)-x; dt=dp(yy); r=dabs(t)  !currant value of
       eps=0.5d0*(dsqrt(dt**2+4.d0*en*r)-dabs(dt))  !autoregularizator
       y=yy-t/(dt+eps)            ! autoregularized Newtonian iterator
             dif=dabs(y-yy)
             dr=dabs(r-rr)
       if(r >= ytol .and. k <= ktol .and. dr > dytol) then
          if(dif <= small) then
             small=dif; ybest=y; rr=r
          endif
          goto 1
       else
          p_newt= ybest
       endif
       drc=drcw
       return
       end
!
!
!
!
!
!
       function p_(x)   ! 2 <= x < infinity 
       implicit real*8(a-h,o-z)
       common/protiarithmi/q(10000000)/wn/wn/ck/vkoch/nn/nn
1      if(x >= 3.d0 .and. x <= q(nn)) then
         if(x <= 1.d3) ixx=floor(x/dlog(x))
         if(x >  1.d3) ixx=floor(dlogintegral(x)-vkoch*dsqrt(x)*dlog(x))
         if(q(ixx) == 0.d0) qq=r(dfloat(ixx))
         if(q(ixx) /= 0.d0) qq=q(ixx)
         if(qq <= x) then
            is=-1; si=-1.d0
         else
            is= 1; si= 1.d0
         endif
       else
         if(x < 3.d0) then
            p_=x-1.d0; return
         endif
       endif
       do i=1, nn
          ii=ixx-is*i
          if(si*(q(ii)-x) <= 0.d0) then
          ix=ii; goto 3; endif
       enddo
3      if(dabs(q(ix)-x) >= dabs(q(ix+is)-x)) ix=ix+is
       xx=q(ix); xi=dfloat(ix); p_= y(ix,xi,x,xx)
       if(vkoch /= 0.d0) goto 1
       contains
       function y(ix,xi,x,xx)   ! left-right inverse quadric pair
        implicit real*8(a-h,o-z)
        common/protiarithmi/q(10000000)/ck/vkoch
              vkoch=0.d0
        if(x <= xx) then
           a=q(ix)-q(ix-1)-1.d0; b=8.d0*a*(q(ix)-x)+1.d0
           if(b <= 0.d0) then
              vkoch=1.d0
           else
              if(a == 0.d0) then
                 y=x-1.d0; else; y=xi+(1.d0-dsqrt(b))/(4.d0*a); endif
           endif
        else
           a=q(ix+1)-q(ix)-1.d0; b=8.d0*a*(x-q(ix))+1.d0
           if(b < 0.d0) then
              vkoch=1.d0
           else
              y=xi+(dsqrt(b)-1.d0)/(4.d0*a)
           endif
        endif
       end function y
       end function p_
!
       function dp_(x)    ! 2 <= x < infinity; inverse of dp(x) 
       implicit real*8(a-h,o-z)
       common/protiarithmi/q(10000000)/wn/wn/ck/vkoch/nn/nn
1      if(x >= 3.d0 .and. x <= q(nn)) then
         if(x <= 1.d3) ixx=floor(x/dlog(x))
         if(x >  1.d3) ixx=floor(dlogintegral(x)-vkoch*dsqrt(x)*dlog(x))
         if(q(ixx) == 0.d0) qq=r(dfloat(ixx))
         if(q(ixx) /= 0.d0) qq=q(ixx)
          if(qq <= x) then
             is=-1; si=-1.d0
          else
             is= 1; si= 1.d0
          endif
       else
         if(x < 3.d0) then
            dp_=1.d0
            return
         endif
       endif
       do i=1, nn
          ii=ixx-is*i
          if(si*(q(ii)-x) <= 0.d0) then
          ix=ii; goto 3; endif
       enddo
3      if(dabs(q(ix)-x) >= dabs(q(ix+is)-x)) ix=ix+is
       xx=q(ix); xi=dfloat(ix); dp_= dy(ix,x,xx)
       if(vkoch /= 0.d0) goto 1
       contains
        function dy(ix,x,xx)   ! left-right inverse quadric pair
        common/protiarithmi/q(10000000)/ck/vkoch
        vkoch=0.d0
        if(x <= xx) then
           a=q(ix)-q(ix-1)-1.d0; b=1.d0+8.d0*a*(q(ix)-x)
           if(b <= 0.d0) then
              vkoch=1.d0
           else
              dy=1.d0/dsqrt(b)
           endif
        else
           a=q(ix+1)-q(ix)-1.d0; b=1.d0+8.d0*a*(x-q(ix))
           if(b <= 0.d0) then
              vkoch=1.d0
           else
              dy=1.d0/dsqrt(b)
           endif
        endif
       end function dy
       end function dp_
!
!
         function r(t)  ! asymptote of function p(t)
!   M. Cipolla, "La determinazione assintotica dell nimo numero primo",
!      Rend. Acad. Sci. Fis. Mat. Napoli, Ser. 3, 8 (1902), 132-166
         implicit real*8(a-h,o-z)
         common/rcorr/rc,drc
       r=t*(dlog(t)+dlog(dlog(t))+(dlog(dlog(t))-2.d0)/              &
         dlog(t)-((dlog(dlog(t)))**2-6.d0*dlog(dlog(t))+11.d0)/      &
         (2.d0*dlog(t)**2)-1.d0)+rc
       end function r

         function dr(x)  ! derivative of asymptote r(t)
         implicit real*8(a-h,o-z)
         common/rcorr/rc,drc
      t1 = dlog(x); t2 = dlog(t1); t3 = t1**2; t4 = t3*t1
      t7 = t1-2.d0; t8 = dlog(t7); t18 = t2**2; t21 = t3**2
      t35 = -56.d0-4.d0*t2*t4+4.d0*t8*t1+50.d0*t1-13.d0*t3+          &
             2.d0*t4-6.d0*t8*t3-26.d0*t2*t1+4.d0*t18*t1+             &
             2.d0*t2*t21+2.d0*t8*t4-1.d0*t18*t3+6.d0*t2*t3+          &
             28.d0*t2-4.d0*t18-4.d0*t21+2.d0*t21*t1
      dr = 0.5d0*t35/t4/t7+drc
      end function dr
!
       subroutine eratosthenes(k)
       implicit real*8(a-h,o-z)
       common/protiarithmi/q(10000000)/rcorr/rc,drc/wn/wn/nn/nn
       common/ck/vkoch
       rc=0.d0; drc=0.d0; np=0; vkoch=0.d0
       do n=2,k
          id=1
          isqrtn=dsqrt(dfloat(n))  
1         id=id+1
          if(n == 2) goto 2
             if(mod(n,id) == 0) goto 3
                if(id >= isqrtn) goto 2
                   goto 1
2               np=np+1; q(np)=n
3            continue
        enddo
        wn=dfloat(np); nn=np; rc=p(wn)-r(wn); drc=dp(wn)-dr(wn)
        end subroutine eratosthenes

       subroutine primes(n)
       implicit real*8(a-h,o-z)
       common/protiarithmi/q(10000000)/rcorr/rc,drc/wn/wn/nn/nn
       common/ck/vkoch
       rc=0.d0; drc=0.d0; nn=n; wn=dfloat(n); vkoch=0.d0
       open(222, file='primes')
        do i=1, n/6+1
           i1=6*i-5;i2=6*i-4;i3=6*i-3;i4=6*i-2;i5=6*i-1;i6=6*i;
           read(222, *) q(i1),q(i2),q(i3),q(i4),q(i5),q(i6)
        enddo
        nn=n; rc=p(wn)-r(wn); drc=dp(wn)-dr(wn)
        end subroutine primes
!
      function dlogintegral(x)
      implicit real*8(a-h,o-z)
      parameter(itol=100, small=1.d-10)
      dli=0.d0; i=0
1     i=i+1; dliw=dli; dm=(dlog(x))**i
      if(dm > small) then
         dli=dli+dfloat(factorial(i-1))/dm
         if(dabs(dliw-dli) > small.and.i <= itol) goto 1
      endif
      dlogintegral=x*dli
      contains
      integer recursive function factorial(l) result(lf)
      integer (4) l
      lf=1
      if(l > 0) then
         lf=l*factorial(l-1); return
      endif
      end function factorial
      end function dlogintegral

       program test_pp_f
       implicit real*8(a-h,o-z)
       open(21, file='p.txt'); open(22, file='p_.txt')
       open(23, file='dp.txt'); open(24, file='dp_.txt')
                         n=6000
       call eratosthenes(n)  ! call primes(n)
!      s=5970.d0; sm=0.05d0; ns=1000
       !uses function primes(n) at n=6000;         returns tab for fig 4
!      s=428.d0; sm=0.01d0; ns=600    ! returns tabs for figs 1, 2 and 3
!      s=154.78d0; sm=0.01d0; ns=1350 !returns tabs for figs 6,7,8 and 9
       s=1.d0; sm=0.1d0; ns=250;     !returns tab for fig 5 
                      ss=s
       do i=1, ns
                      s=s+sm
          write(21,*) s,  p(s)
          write(23,*) s, dp(s)
       enddo
                      sm=(p(s)-p(ss))/dfloat(ns); s=p(ss)
       do i=1, ns
                      s=s+sm
          write(22,*) s, p_(s)       ! p_newt(s)
          write(24,*) s, dp_(s)
       enddo
       end


\end{verbatim}

\large
{\bf Application 2.}\\
 Solution of the equation $x(1)^2+x(2)^2=x(3)^2+1$ over primes
\small
\begin{verbatim}
!---user---module---to---the---main---program---afxy--------------------
      subroutine FXY(m,n,np,neq,f,x,pp,df,yr)
      implicit real*8(a-h,o-z)
      DIMENSION X(1),pp(1),DF(1),YR(1)
      COMMON/LSMH/LSMH,MQH,NSOLH/BXH/D1,D2,BL(600),BR(600)
      COMMON/RETFH/LF1,LF2,LF3,NDAT/FR/FR/SLMH/SSVH,S3H
           go to (1,2), np
!-----run-time--section-------------------------------------------------
   2  pi=dacos(-1.D0)
      goto(21,22), neq
!.first..equation.......................................................
  21  f=x(1)**2+x(2)**2-x(3)**2-1.d0              ! diophantine equation
      df(1)= 2.d0*x(1); df(2)= 2.d0*x(2); df(3)=-2.d0*x(3) ! derivatives
      return
!.second..equation......................................................
  22  f=0.d0
      do i=1,3
         df(i)=pi*dsin(2.d0*pi*p_(x(i)))*dp_(x(i))         ! derivatives
         f=f+(dsin(pi*p_(x(i))))**2  ! from-reals-to-primes equation
      enddo
      return
!-------pre-execution--section--set--the--controls--to--afxy-program---
   1  nn=10000; call eratosthenes(nn)    ! prime number table creation
           n=3                                    ! number of unknowns 
           m=2                                    ! number of equations
                np=-3    ! produce autoregularized Gauss-Newton process
                f=1.D-34    ! accuracy level for the residual f(x)-y
                lsmh=1                         ! Gene Golub's SVD-method
                  ssvh=1.d-16 ! minimal characterist numb in SVD-method
                mqh=1                  ! Jorge More's addaptive scaling
                nsolh=20               ! limit of the sought solutions
           do i=1, n
              x(i)=20.d0               ! guesses
              bl(i)=2.d0; br(i)=100.d0 ! constraints (definition domain Df)
           enddo
          lf1=1; lf2=2
          return
      end
!
          FOUND SOLUTIONS:
 Solution #  1 (Plan  2; x0 #  2; eps0=2.88D-08) :
 K=  32-----------------------------------------------------------------------
 (Df)`(fx-y)    fx-y      chi.sqr.   (Df)`Df    l.p.prec.     eps      en/r
 2.66223D-21 6.39687D-15 4.09199D-29 1.3172D+04 3.5976D-23 0.000D+00 4.8D-06
 Unknowns:
 X(  1)= 2.8999999979D+01 X(  2)= 2.3000000010D+01 X(  3)= 3.6999999989D+01

 Solution #  2 (Plan  2; x0 #  2; eps0=2.88D-03) :
 K=  26-----------------------------------------------------------------------
 (Df)`(fx-y)    fx-y      chi.sqr.   (Df)`Df    l.p.prec.     eps      en/r
 5.32001D-18 9.98342D-13 9.96687D-25 1.6260D+04 1.2170D-13 0.000D+00 4.8D-02
 Unknowns:
 X(  1)= 2.9000000270D+01 X(  2)= 2.8999999962D+01 X(  3)= 4.1000000164D+01

 Solution #  3 (Plan  2; x0 #  2; eps0=2.88D-01) :
 K=  50-----------------------------------------------------------------------
 (Df)`(fx-y)    fx-y      chi.sqr.   (Df)`Df    l.p.prec.     eps      en/r
 7.48702D-21 1.30119D-14 1.69311D-28 2.6924D+04 7.0632D-23 0.000D+00 5.3D-01
 Unknowns:
 X(  1)= 4.2999999971D+01 X(  2)= 3.1000000004D+01 X(  3)= 5.2999999979D+01

 Solution #  4 (Plan  2; x0 #  3; eps0=2.42D-08) :
 K=  35-----------------------------------------------------------------------
 (Df)`(fx-y)    fx-y      chi.sqr.   (Df)`Df    l.p.prec.     eps      en/r
 1.74942D-20 2.45302D-14 6.01732D-28 1.3197D+04 1.3640D-13 0.000D+00 3.1D-05
 Unknowns:
 X(  1)= 2.2999999992D+01 X(  2)= 2.8999999964D+01 X(  3)= 3.6999999967D+01

 Solution #  5 (Plan  2; x0 #  3; eps0=2.42D-07) :
 K=  52-----------------------------------------------------------------------
 (Df)`(fx-y)    fx-y      chi.sqr.   (Df)`Df    l.p.prec.     eps      en/r
 9.75966D-22 3.70459D-15 1.37240D-29 5.0600D+03 1.7760D-14 0.000D+00 3.1D-04
 Unknowns:
 X(  1)= 1.9000000009D+01 X(  2)= 1.3000000011D+01 X(  3)= 2.3000000013D+01

 Solution #  6 (Plan  2; x0 #  3; eps0=2.42D-04) :
 K=  32-----------------------------------------------------------------------
 (Df)`(fx-y)    fx-y      chi.sqr.   (Df)`Df    l.p.prec.     eps      en/r
 6.32452D-21 1.18087D-14 1.39445D-28 2.1246D+04 6.7279D-23 0.000D+00 3.1D+00
 Unknowns:
 X(  1)= 2.9000000000D+01 X(  2)= 3.6999999973D+01 X(  3)= 4.6999999979D+01

 Solution #  7 (Plan  2; x0 #  4; eps0=2.01D-06) :
 K=  52-----------------------------------------------------------------------
 (Df)`(fx-y)    fx-y      chi.sqr.   (Df)`Df    l.p.prec.     eps      en/r
 9.97306D-21 1.59910D-14 2.55714D-28 2.0899D+04 1.0602D-22 0.000D+00 3.1D-05
 Unknowns:
 X(  1)= 2.3000000008D+01 X(  2)= 4.0999999968D+01 X(  3)= 4.6999999976D+01

 Solution #  8 (Plan  2; x0 #  4; eps0=2.01D-05) :
 K=  56-----------------------------------------------------------------------
 (Df)`(fx-y)    fx-y      chi.sqr.   (Df)`Df    l.p.prec.     eps      en/r
 2.27957D-20 3.09817D-14 9.59867D-28 1.6685D+04 2.5783D-13 0.000D+00 3.1D-03
 Unknowns:
 X(  1)= 1.3000000029D+01 X(  2)= 4.1000000030D+01 X(  3)= 4.3000000037D+01
!
!
 Solution #  9 (Plan  2; x0 #  5; eps0=1.74D-08) :
 K=  40-----------------------------------------------------------------------
 (Df)`(fx-y)    fx-y      chi.sqr.   (Df)`Df    l.p.prec.     eps      en/r
 2.16010D-16 1.23358D-11 1.52172D-22 2.6935D+04 1.5541D-13 0.000D+00 1.1D-06
 Unknowns:
 X(  1)= 3.0999999926D+01 X(  2)= 4.3000000887D+01 X(  3)= 5.3000000676D+01

 Solution # 10 (Plan  2; x0 #  5; eps0=1.74D-04) :
 K=  37-----------------------------------------------------------------------
 (Df)`(fx-y)    fx-y      chi.sqr.   (Df)`Df    l.p.prec.     eps      en/r
 4.52579D-21 1.02193D-14 1.04433D-28 8.8041D+03 2.8456D-14 0.000D+00 1.1D-02
 Unknowns:
 X(  1)= 2.9000000019D+01 X(  2)= 1.1000000013D+01 X(  3)= 3.1000000022D+01

 Solution # 11 (Plan  2; x0 #  5; eps0=1.74D+00) :
 K=  29-----------------------------------------------------------------------
 (Df)`(fx-y)    fx-y      chi.sqr.   (Df)`Df    l.p.prec.     eps      en/r
 4.96913D-21 1.02836D-14 1.05752D-28 2.0492D+04 2.3436D-13 0.000D+00 2.2D+01
 Unknowns:
 X(  1)= 4.3000000024D+01 X(  2)= 1.8999999989D+01 X(  3)= 4.7000000018D+01

 Solution # 12 (Plan  2; x0 #  6; eps0=2.43D-03) :
 K=  44-----------------------------------------------------------------------
 (Df)`(fx-y)    fx-y      chi.sqr.   (Df)`Df    l.p.prec.     eps      en/r
 1.26391D-23 1.87462D-16 3.51422D-32 5.0621D+03 3.5534D-14 0.000D+00 1.3D+01
 Unknowns:
 X(  1)= 1.3000000000D+01 X(  2)= 1.9000000003D+01 X(  3)= 2.3000000003D+01

 Solution # 13 (Plan  2; x0 #  7; eps0=2.70D-08) :
 K=  58-----------------------------------------------------------------------
 (Df)`(fx-y)    fx-y      chi.sqr.   (Df)`Df    l.p.prec.     eps      en/r
 7.20035D-21 1.26460D-14 1.59923D-28 4.3148D+04 4.7940D-16 0.000D+00 6.7D-03
 Unknowns:
 X(  1)= 4.1000000029D+01 X(  2)= 5.2999999979D+01 X(  3)= 6.7000000001D+01

 Solution # 14 (Plan  2; x0 #  7; eps0=2.70D-02) :
 K=  44-----------------------------------------------------------------------
 (Df)`(fx-y)    fx-y      chi.sqr.   (Df)`Df    l.p.prec.     eps      en/r
 7.26781D-17 6.74038D-12 4.54328D-23 6.5404D+04 2.1514D-13 0.000D+00 6.7D-01
 Unknowns:
 X(  1)= 7.1000000315D+01 X(  2)= 4.3000000534D+01 X(  3)= 8.3000000546D+01

 Solution # 15 (Plan  2; x0 #  9; eps0=2.73D-03) :
 K=  28-----------------------------------------------------------------------
 (Df)`(fx-y)    fx-y      chi.sqr.   (Df)`Df    l.p.prec.     eps      en/r
 2.31418D-21 6.77216D-15 4.58621D-29 8.8041D+03 3.1949D-14 0.000D+00 1.6D+00
 Unknowns:
 X(  1)= 1.1000000015D+01 X(  2)= 2.9000000013D+01 X(  3)= 3.1000000017D+01
!
!
 Solution # 16 (Plan  2; x0 # 10; eps0=2.78D-05) :
 K=  71-----------------------------------------------------------------------
 (Df)`(fx-y)    fx-y      chi.sqr.   (Df)`Df    l.p.prec.     eps      en/r
 1.32835D-20 2.01058D-14 4.04243D-28 4.3148D+04 9.9131D-23 0.000D+00 1.7D-04
 Unknowns:
 X(  1)= 5.2999999967D+01 X(  2)= 4.0999999995D+01 X(  3)= 6.6999999970D+01

 Solution # 17 (Plan  2; x0 # 13; eps0=2.58D+00) :
 K=  36-----------------------------------------------------------------------
 (Df)`(fx-y)    fx-y      chi.sqr.   (Df)`Df    l.p.prec.     eps      en/r
 7.87327D-20 5.71405D-14 3.26504D-27 7.4410D+04 2.8416D-13 0.000D+00 8.1D+01
 Unknowns:
 X(  1)= 7.9000000002D+01 X(  2)= 4.0999999930D+01 X(  3)= 8.8999999970D+01

 Solution # 18 (Plan  2; x0 # 15; eps0=2.47D-03) :
 K=  55-----------------------------------------------------------------------
 (Df)`(fx-y)    fx-y      chi.sqr.   (Df)`Df    l.p.prec.     eps      en/r
 1.70828D-20 2.50831D-14 6.29162D-28 4.7012D+04 5.1160D-13 0.000D+00 3.4D+00
 Unknowns:
 X(  1)= 7.1000000031D+01 X(  2)= 1.7000000020D+01 X(  3)= 7.3000000035D+01

 Solution # 19 (Plan  2; x0 # 16; eps0=2.29D-09) :
 K=  47-----------------------------------------------------------------------
 (Df)`(fx-y)    fx-y      chi.sqr.   (Df)`Df    l.p.prec.     eps      en/r
 6.22186D-21 1.25141D-14 1.56603D-28 1.6703D+04 1.2445D-13 0.000D+00 1.4D-05
 Unknowns:
 X(  1)= 4.1000000024D+01 X(  2)= 1.3000000008D+01 X(  3)= 4.3000000025D+01

 Solution # 20 (Plan  2; x0 # 21; eps0=1.97D-08) :
 K=  34-----------------------------------------------------------------------
 (Df)`(fx-y)    fx-y      chi.sqr.   (Df)`Df    l.p.prec.     eps      en/r
 1.02797D-20 1.75525D-14 3.08090D-28 2.1356D+04 1.5072D-13 0.000D+00 2.1D-05
 Unknowns:
 X(  1)= 4.0999999971D+01 X(  2)= 2.2999999991D+01 X(  3)= 4.6999999971D+01

     Singular values of the matrix Z=(Df)`Df at the best approximation X( 34)
     (errSVD=  0):

     1) 0.0000000000000D+00   2) 2.0194839173658D-28   3) 1.8088976985834D+04

     Cond(Z(x( 34))=104 ;  GMcond(Z(x( 34)))= 45 ;  nullity(Z(x( 34)))=  2 .

 Residual_tab of the vector r=y-f(x( 34)) (M=   2):
 0.0D+00,  1;  -1.8D-14,  1;

   Final report :   202 iterative processes were tried,
                   9093 times subroutine FXY was called,
                   8911 iterations were produced, and
                     20 solutions were found.
-------------------------------------------------------------------------------
\end{verbatim}

\newpage

\large
{\bf Application 3.}\\
 Solution of the system {$x(1)^2+x(2)^2=x(3)^2+1, x(3)-x(1)=2$}
 over twins
\small
\begin{verbatim}
!---user's---module--to--the--main--program--afxy---------------------------
      subroutine FXY(m,n,np,neq,f,x,pp,df,yr)
      implicit real*8(a-h,o-z)
      DIMENSION X(1),pp(1),DF(1),YR(1)
      COMMON/LSMH/LSMH,MQH,NSOLH/BXH/D1,D2,BL(600),BR(600)
      COMMON/RETFH/LF1,LF2,LF3,NDAT/SLMH/SSVH,S3H
           go to (1,2), np
!-----run-time--section-------------------------------------------------
   2       continue; pi=dacos(-1.D0)
      goto(21,22,23), neq
!...first..equation.....................................................
  21  f=x(1)**2+x(2)**2-x(3)**2-1.d0              ! diophantine equation
      df(1)= 2.d0*x(1); df(2)= 2.d0*x(2); df(3)=-2.d0*x(3) ! derivatives
      return
!...second..equation....................................................
  22  f=x(3)-x(1)-2.d0                              ! twin-pair equation
      df(1)=-1.d0; df(2)= 0.d0; df(3)= 1.d0                ! derivatives
      return
!...third..equation.....................................................
  23  f=0.d0
      do i=1,3
         df(i)=pi*dsin(2.d0*pi*p_(x(i)))*dp_(x(i))         ! derivatives
         f=f+(dsin(pi*p_(x(i))))**2  ! from-reals-to-primes equation
      enddo
      return
!-------pre-execution--section--set--the--controls--to--afxy-program---
   1       continue
           nn=10000; call eratosthenes(nn)! prime number table creation
           n=3                                    ! number of unknowns 
           m=3                                    ! number of equations
                np=-3    ! produce autoregularized Gauss-Newton process
                f=1.D-34      ! accuracy level for the residual f(x)-y
                  lsmh=1                       ! Gene Golub's SVD-method
                ssvh=1.d-16   ! minimal characterist numb in SVD-method
                mqh=1                  ! Jorge More's addaptive scaling
                nsolh=10               ! limit of the sought solutions
           do i=1, n
              x(i)=50.d0               ! guesses
              bl(i)=2.d0; br(i)=101.d0 ! constraints (definition domain Df)
           enddo
          lf1=1; lf2=2
          return
      end
!
!
!
!
!
!
       FOUND SOLUTIONS:
 Solution #  1 (Plan  2; x0 #  1; eps0=3.00D-02) :
 K=  31-----------------------------------------------------------------------
 (Df)`(fx-y)    fx-y      chi.sqr.   (Df)`Df    l.p.prec.     eps      en/r
 9.94498D-20 7.96082D-14 6.33746D-27 4.7014D+04 6.3947D-14 0.000D+00 3.6D-01
 Unknowns:
 X(  1)= 7.0999999937D+01 X(  2)= 1.6999999993D+01 X(  3)= 7.2999999937D+01

 Solution #  2 (Plan  2; x0 #  1; eps0=3.00D+00) :
 K=  36-----------------------------------------------------------------------
 (Df)`(fx-y)    fx-y      chi.sqr.   (Df)`Df    l.p.prec.     eps      en/r
 1.84559D-21 5.58920D-15 3.12392D-29 1.6686D+04 7.1257D-15 0.000D+00 7.2D+00
 Unknowns:
 X(  1)= 4.1000000017D+01 X(  2)= 1.3000000003D+01 X(  3)= 4.3000000017D+01

 Solution #  3 (Plan  2; x0 #  2; eps0=2.88D-08) :
 K=  34-----------------------------------------------------------------------
 (Df)`(fx-y)    fx-y      chi.sqr.   (Df)`Df    l.p.prec.     eps      en/r
 9.06718D-24 1.65074D-16 2.72495D-32 4.7800D+02 1.7750D-15 0.000D+00 3.1D-04
 Unknowns:
 X(  1)= 5.0000000028D+00 X(  2)= 5.0000000011D+00 X(  3)= 7.0000000028D+00

 Solution #  4 (Plan  2; x0 #  2; eps0=2.88D-07) :
 K=  34-----------------------------------------------------------------------
 (Df)`(fx-y)    fx-y      chi.sqr.   (Df)`Df    l.p.prec.     eps      en/r
 1.14803D-22 8.85711D-16 7.84485D-31 1.6145D+03 3.5507D-15 0.000D+00 3.1D-04
 Unknowns:
 X(  1)= 1.0999999993D+01 X(  2)= 6.9999999981D+00 X(  3)= 1.2999999993D+01

 Solution #  5 (Plan  2; x0 #  2; eps0=2.88D-03) :
 K=  33-----------------------------------------------------------------------
 (Df)`(fx-y)    fx-y      chi.sqr.   (Df)`Df    l.p.prec.     eps      en/r
 4.85228D-22 2.29735D-15 5.27781D-30 8.8060D+03 3.5523D-15 0.000D+00 3.1D-01
 Unknowns:
 X(  1)= 2.9000000011D+01 X(  2)= 1.1000000002D+01 X(  3)= 3.1000000011D+01

   Final report :  300 iterative processes were tried,
                  6108 times subroutine FXY was called,
                  6027 iterations were produced, and
                     5 solutions were found.
-------------------------------------------------------------------------------
\end{verbatim}

\end{document}